\documentclass[10pt]{article}

\newcommand{\de}{\delta}
\newcommand{\la}{\lambda}

\newcommand{\ra}{\rightarrow}

\DeclareSymbolFont{AMSb}{U}{msb}{m}{n}
\DeclareMathSymbol{\N}{\mathbin}{AMSb}{"4E}
\DeclareMathSymbol{\Z}{\mathbin}{AMSb}{"5A}
\DeclareMathSymbol{\R}{\mathbin}{AMSb}{"52}
\DeclareMathSymbol{\Q}{\mathbin}{AMSb}{"51}
\DeclareMathSymbol{\I}{\mathbin}{AMSb}{"49}
\DeclareMathSymbol{\C}{\mathbin}{AMSb}{"43}

\usepackage{graphicx}
\usepackage{amssymb}
\usepackage{epstopdf}
\usepackage{amsmath}
\usepackage{amsfonts}

\newtheorem{theorem}{Theorem}

\newtheorem{proposition}{Proposition}
\newtheorem{lemma}{Lemma}

\def\phi{{\varphi}}

\begin{document}

\title{On the existence of Hamiltonian paths connecting Lagrangian submanifolds}
\author{ Nassif  Ghoussoub\thanks{Partially supported by a grant
from the Natural Sciences and Engineering Research Council of Canada.  } \quad  and \quad Abbas Moameni\thanks{Research supported by a postdoctoral fellowship at the University of British Columbia.}
\\
\small Department of Mathematics,
\small University of British Columbia, \\
\small Vancouver BC Canada V6T 1Z2 \\
\small {\tt nassif@math.ubc.ca} \\
\small {\tt moameni@math.ubc.ca}
\\
}
\maketitle

\abstract

We use a new variational method --based on the theory of anti-selfdual Lagrangians developed in \cite{G2} and \cite{G3}-- to establish the existence of solutions of convex Hamiltonian systems that connect two given Lagrangian submanifolds in $\R^{2N}$. We also consider the case where the Hamiltonian is only semi-convex. A variational principle  is also used to establish existence for  the corresponding Cauchy problem. The case of periodic solutions will be considered in a forthcoming paper \cite{GM2}.

\section{Introduction}
We consider the following Hamiltonian System
\begin{eqnarray}
\left\{\begin{array}{rcll}
\dot{p}(t) &\in &\partial_2 H\big( p(t),q(t)\big)\quad &t\in (0,T),\\
-\dot{q} (t) &\in &\partial_1 H\big( p(t),q(t)\big) &t\in (0,T),
\end{array}\right.
\end{eqnarray}
where $H:\R ^N\times\R ^N\ra\R$ is a convex and lower semi
continuous function and $T>0$. We develop a new variational approach to establish existence of solutions  satisfying two types of boundary conditions. The first one requires the path to connect  two Lagrangian submanifolds associated to given convex  lower semi continuous functions $\psi_1$ and $\psi_2$  on $\R^N$, that is
\begin{eqnarray}
q(0)\in\partial\psi_1\big( p(0)\big)\quad {\rm and}\quad
-p(T)\in\partial\psi_2\big( q(T)\big).
\end{eqnarray}
In other words, the Hamiltonian path must connect the graph of $\partial \psi_1$ to the graph of $-\partial \psi_2$.
The second is simply an initial value problem of the form
\begin{eqnarray}
p(0)=p_0,\quad q(0)=q_0
\end{eqnarray}
where $p_0$ and $q_0$ are two given vectors in $\R ^N$.

The solutions will be obtained from a novel variational principle developed  in full generality in a series of papers \cite{G2}, \cite{G3} and \cite{GT1}. It is based on the concept of anti-selfdual Lagrangians to which one associates action functionals whose infimum is necessarily equal to zero. The equations are then derived from the limiting case in Legendre-Fenchel duality as opposed to standard Euler-Lagrange theory.

 In the next section, we start with the case of convex Hamiltonian systems connecting Lagrangian submanifolds. This is then extended to the semi-convex case in section 4.  The corresponding Cauchy problem is studied in section 3.

\section{Connecting Lagrangian submanifolds}

Given a time $T>0$, we let $X=W^{1,2}(0,T;\R ^N)$ be the one-dimensional Sobolev space
endowed with the norm
$\| u\| ={\left(\| u\|_{L^2}^2 +{\| \dot{u}\|}_{L^2}^2
\right)}^{\frac{1}{2}}$
where$\| u\|_{L^2}={\left(\int_0^T |u|^2 dt\right)}^{\frac{1}{2}}$ stands for the norm on $L^2:=L^2(0,T;\R ^N)$. For every $p,q\in\R ^N$, $p\cdot q$ denotes the inner product in $\R
^N$ and $(p,q)\cdot (r,s)$ denotes the inner product in $\R
^N\times\R ^N$ defined by $ (p,q)\cdot (r,s)=p\cdot r+q\cdot s. $
 
Say that a Hamiltonian   $H$ on $\R^{2N}$ is  $\beta$-subquadratic for $\beta >0$, if  for some positive constants $\alpha,  \gamma$, we have,
\begin{eqnarray}
\hbox{$-\alpha \leq H(p,q)\leq \frac {\beta}{2} (|p|^2+
|q|^2)+\gamma$ \quad for all $(p,q)\in \R^{2N}$.}
\end{eqnarray}
We shall prove the following result.
\begin{theorem}
Suppose $H:\R ^{2N}\ra\R$ is a  convex lower
semi-continuous $\beta$-subquadratic Hamiltonian with
\begin{equation}
\label{beta}
 \beta < \frac{1}{2\max(2T^2,1)}.
 \end{equation}
Let $\psi_1$ and $\psi_2$ be two  convex lower
semi-continuous and coercive  functions on $R ^N$  such that one of
them satisfies the following condition:
\begin{eqnarray}\begin{array}{l}
\liminf\limits_{|p|\ra +\infty}\frac{\psi_i (p)}{|p|^2} > 2T \quad {\rm for}\quad i=1 \,\,  \text{\rm or} \,\, 2.\\
\end{array}
\end{eqnarray}
Then the minimum of the  functional
\begin{eqnarray*}
I(p,q):&=&
\int_0^T \big[ H\big( p(t),q(t)\big) + H^*\big(-\dot{q} (t),\dot{p}(t)\big) +2\dot{q}(t)\cdot p(t)\big] \, dt\\
  &&\quad + \psi_2 \big( q(T)\big) +\psi_2^* \big(  -p(T)\big) +\psi_1 \big( p(0)\big) +\psi_1^* \big( q(0))
\end{eqnarray*}
on $Y=X\times X$ is zero and is attained at a solution of 
\begin{eqnarray}
\label{connecting}
\left\{\begin{array}{rcll}
\dot{p}(t) &\in &\partial_2 H\big( p(t),q(t)\big)\quad &t\in (0,T),\\
-\dot{q} (t) &\in &\partial_1 H\big( p(t),q(t)\big) &t\in (0,T),\\
q(0)\in\partial\psi_1\big( p(0)\big)&\&
&
-p(T)\in\partial\psi_2\big( q(T)\big).
\end{array}\right.
\end{eqnarray}
\end{theorem}
\bigskip

Before we proceed with the proof, we note that condition (\ref{beta}) is satisfied as soon as we have 
 \begin{eqnarray*}
-\alpha \leq H(p,q)\leq \beta (|p|^r+ |q|^r+1) \quad \quad
(1< r<2)
\end{eqnarray*}
 where  $\alpha, \beta$ are any positive constants.
 
The proof requires a few preliminary lemmas, but first and anticipating that at some point of the proof,  the conjugate $H^*$ of $H$ needs to be finite everywhere (i.e, $H$ coercive), we start by replacing $H$ with the following perturbed Hamiltonian 
$H_{\epsilon}(p,q)= \frac
{\epsilon}{2}(|p|^2+|q|^2) + H(p,q)$ 
for some $\epsilon >0$. It is then clear that
\begin{eqnarray}
\frac {1}{2(\beta+\epsilon)}(|p|^2+|q|^2)-\gamma \leq
H^*_{\epsilon}(p,q)\leq  \frac {1}{2 \epsilon}(|p|^2+|q|^2)+\alpha .
\end{eqnarray}
 
\begin{lemma} For any convex Hamiltonian $H$, and convex  lower semi-continuous functions $\psi_1$, $\psi_2$,   we have that $I(p,q)\geq 0$ for
every $(p,q)\in X\times X$.
\end{lemma}
\noindent{\bf Proof:} Use that
\begin{eqnarray*}
2\int_0^T\dot{q}\cdot p\, dt=\int_0^T\dot{q}\cdot p\, dt-\int_0^T
q\cdot\dot{p}\, dt+q(T)\cdot p(T)-q(0)\cdot p(0)
\end{eqnarray*}
to write
\begin{eqnarray*}
I(p,q):&=&
\int_0^T \big[ H\big( p(t),q(t)\big) + H^*\big( -\dot{q} (t),\dot{p}(t)\big) +2\dot{q}(t)\cdot p(t)\big] \, dt\\
  &&\quad + \psi_2 \big( p(T)\big) +\psi_2^* \big(  -q(T)\big) +\psi_1 \big( p(0)\big) +\psi_1^* \big( q(0))\\
  &=&
\int_0^T \big[ H\big( p(t),q(t)\big) + H^*\big( -\dot{q} (t),\dot{p}(t)\big) +\dot{q}(t)\cdot p(t)-\dot{p}(t)\cdot q(t)\big] \, dt\\
  &&\quad + \big[\psi_2 \big( q(T)\big) +\psi_2^* \big(  -p(T)\big)+p(T)\cdot q(T)
  \big]\\
  &&\quad + \big[\psi_1 \big( p(0)\big) +\psi_1^* \big( q(0)\big) -p(0)\cdot q(0) \big]\\
  &\geq&0
\end{eqnarray*}
by applying three times the Legendre-Fenchel inequality.

\begin{lemma} To a convex lower semi-continuous Hamiltonian $H$ on $\R^{2N}$, we associate the following ``Functional Lagrangian"  $L:Y\times Y\ra \R\cup\{
+\infty\}$ defined as
\begin{eqnarray*}
L(r,s;p,q):&=& \int_0^T\left[ (\dot{r},-\dot{s})\cdot (q,p)+H^* (-\dot{q},\dot{p})-H^*(-\dot{s},\dot{r})+2\dot{q}\cdot p\right]\, dt\\
&&- p(T)\cdot s(T)+\psi_2\big( q(T)\big) -\psi_2 \big( s(T)\big)
+ r(0)\cdot q(0)+\psi_1\big( p(0)\big) -\psi_1 \big( r(0)\big).
\end{eqnarray*}
Then, we have $I({p},{q})\leq \sup\limits_{(r,s)\in X\times
X}L(r,s;{p},{q})$ for every $({p},{q}) \in Y$. 
\end{lemma}

 \noindent{\bf Proof:} Indeed, set
\begin{eqnarray*}
A:=\left\{ r\in X:r(t)=\int_0^t f(\alpha)\, d\alpha+y,\mbox{ for some $y\in\R
^N$ and }f\in L_2(0,T;\R ^N)\right\}
\end{eqnarray*}
and
\begin{eqnarray*}
B:=\left\{ s\in X:s(t)=-\int_t^T g(\alpha)\, d\alpha+x,\mbox{ for
some $x\in\R ^N$ and }g\in L_2(0,T;\R ^N)\right\}
\end{eqnarray*}
and note that for every $({p},{q}) \in X\times X$ we have
\begin{eqnarray*}
\sup\limits_{X\times
X}L(r,s;{p},{q})
&\geq& \sup_{(r,s)\in A\times B}L(r,s;{p},{q})\\
  &=& \sup_{\stackrel{f,g\in L_2}{x,y\in\R ^N}}
   \big\{
   \int_0^T \left[ (f,-g)\cdot ({q},{p})+H^*(-\dot{q},\dot{p})
    -H^*(-g,f)+2\dot{{q}}\cdot {p}\right]\, dt\big\}.\\
  &&\quad\quad \quad \quad \quad -x\cdot{p}(T)+\psi_2\big({q}(T)\big) -\psi_2 (x)+ y\cdot{q} (0)+\psi_1\big({p}(0)\big) -\psi_1(y)\big\} \\
  &=& \sup_{f,g\in L^2}\left\{\int_0^T \left[(f,-g)\cdot ({q},{p})+H^*(-\dot{{q}},\dot{{p}})-H^*(-g,f)+2\dot{{q}}\cdot{p}\right]\, dt\right\}\\
  && \quad +\sup_{x\in\R ^N}\left\{ -x\cdot{p}(T)+\psi_2\big({q}(T)\big) -\psi_2(x)\right\} +\sup_{y\in\R ^N}\left\{ y\cdot{q}(0)+\psi_1\big({p}(0)\big) -\psi_1(y)\right\}\\
  &=& \int_0^T\left[ H^*(-\dot{{q}},\dot{{p}})+H({p},{q})+2\dot{{q}}\cdot{p}\right]\, dt\\
  &&\quad + \left\{ \psi_2\big({q}(T)\big) +\psi_2^*\big( -{p}(T)\big)\right\} +\left\{\psi_1\big({p}(0)\big) +\psi_1^*\big({q}(0)\big)\right\}\\
  &=& \int_0^T\left[ H^*(-\dot{{q}},\dot{{p}})+H({p},{q})+\dot{{q}}\cdot{p}-{q}\cdot \dot{{p}}\right]\, dt\\
 && \quad + \big[\psi_2\big({q}(T)\big) +\psi_2^*\big( -{p}(T)\big)+{p}(T)\cdot{q}(T)\big]\\
 &&\quad + \big[ \psi_1\big({p}(0)\big) +\psi_1^*\big({q}(0)\big) -{p}(0)\cdot
  {q}(0)\big]\\
  &=& I ({p},{q}).
\end{eqnarray*}
 
\begin{lemma} Under the above conditions, and assuming that $\epsilon$ is small enough so that $\beta +\epsilon< \frac{1}{2\max(2T^2,1)}$, then we have the following
coercivity property:
\[
\hbox{$L_\epsilon(0,0;p,q)\ra +\infty$ as $
 \| p\| +\| q\|\ra +\infty$.}
\]
where $L_\epsilon$ is the functional Lagrangian associated to the perturbed Hamiltonian $H_\epsilon$. 
\end{lemma}
\noindent{\bf Proof:}  Without loss of generality we assume $\psi_1$
satisfies (6). An easy calculation shows that
\begin{eqnarray}
\| p\|_{L^2}\leq T\| \dot{p}\|_{L^2} +\sqrt{T} \big( |p(0)|, \quad
\| q\|_{L^2}\leq T\|\dot{q}\|_{L^2}+\sqrt{T} |q(T)|.
\end{eqnarray}
Also note that  $\frac {1}{2(\beta+\epsilon)}(|p|^2+|q|^2)-\gamma
\leq H^*_{\epsilon}(p,q),$ hence modulo a constant we have
\begin{eqnarray}
L(0,0,p,q)\geq \frac {1}{2(\beta+\epsilon)}\int_0^T (|\dot{q}|^2 +
|\dot{p}|^2)\, dt+2\int_0^T \dot{q}\cdot p\, dt+\psi_2\big(
q(T)\big) +\psi_1\big( p(0)\big)
\end{eqnarray}
Using Holder's inequality and inequality  (9) for the second term on
the right hand side of (10), imply
\begin{eqnarray}
\left|\int_0^T\dot{q}\cdot p\, dt\right| \leq  \frac{1}{2}\int |p|^2\, dt +\frac{1}{2}\int |\dot{q}|^2\, dt
\leq  T^2\int |\dot{p}|^2\, dt +T|p(0)|^2+\frac{1}{2}\int
|\dot{q}|^2\, dt.
\end{eqnarray}
  From (10) and (11), we get
\begin{eqnarray*}
L(0,0,p,q) &\geq & \frac {1}{2(\beta+\epsilon)}\int_0^T
(|\dot{q}|^2+  |\dot{p}|^2)\, dt-\max \left( 2T^2,1\right)\int_0^T
(|\dot{q}|^2+|\dot{p}|^2)dt
\\
&& + \psi_2\big( q(T)\big) +\psi_1 \big( p(0)\big) -2T |p(0)|^2
\end{eqnarray*}
which together with the coercivity condition on $\psi_1$ and $\psi_2$
 and the fact that $\beta+ \epsilon < \frac{1}{2 \max(2T^2,1)}$ imply the claimed coercivity for $L$.\\

\noindent The theorem is now a consequence of the following Ky-Fan type min-max theorem which is essentially due to Brezis-Nirenberg-Stampachia (see \cite{BNS}).
 \begin{lemma} \label{KF} Let $D$ be an unbounded  closed convex subset of  a reflexive Banach space $Y$,   and let $L(x,y)$ be a real valued function  on $D\times D $ that satisfies the following conditions:
\begin{description}
\item (1) $L(x,x) \leq 0$ for every $x\in D$.
\item (2) For each $x\in D$,  the function  $y \to L(x,y)$ is concave.
\item (3) For each $y\in D$, the function $x\to L(x,y)$ is weakly lower semi-continuous.
\item (4) The set $D_0=\{x\in D; L(x,0)\leq 0\}$ is bounded in $Y$.
 \end{description}
Then there exists $x_{0}\in D$ such that $\sup\limits_{y\in D}L(x_{0},y)\leq 0$.
\end{lemma}

\noindent {\bf Proof of Theorem 1:} It is easy to see that $\tilde L_\epsilon$ defined by
\begin{eqnarray*}
\tilde L_\epsilon (p,q; r,s)&=&L_\epsilon (r,s;p,q)\\
&=&\int_0^T\left[ (\dot{r},-\dot{s})\cdot (q,p)+H_\epsilon^* (-\dot{q},\dot{p})-H_\epsilon^*(-\dot{s},\dot{r})+2\dot{q}\cdot p\right]\, dt\\
&&- p(T)\cdot s(T)+\psi_2\big( q(T)\big) -\psi_2 \big( s(T)\big)
+ r(0)\cdot q(0)+\psi_1\big( p(0)\big) -\psi_1 \big( r(0)\big).
\end{eqnarray*}
satisfies  all the hypothesis of lemma 4 on the space $Y=X\times X$. Indeed, From (8) $L_\epsilon$ is real valued and it is clear that $L_\epsilon (p,q; p,q)=0$ and Lemma 3 gives that the set $Y_0=\{ (p,q)\in Y; L_\epsilon (0,0;p,q)\leq 0\}$ is bounded in $Y$. The function $(r,s) \to L_\epsilon (r,s; p,q)$ is concave for every $(p,q)$ while $(p,q) \to L_\epsilon (r,s; p,q)$ is weakly lower semi-continuous for every $(r,s)\in Y$. \\
It follows that there exists $({ p_{\epsilon}}, { q_{\epsilon}})$
such that $\sup\limits_{(r,s)\in X\times X}L_\epsilon (r,s; { p_{\epsilon}}, {
q_{\epsilon}})\leq 0$, so that by Lemma (2) we have
\[
 I({ p_{\epsilon}},{ q_{\epsilon}}) \leq \sup\limits_{(r,s)\in X\times X}L_\epsilon (r,s; { p_{\epsilon}}, { q_{\epsilon}}) \leq 0.
\]
 On the other hand by Lemma (1) we have that $ I_\epsilon ({ p_{\epsilon}},{ q_{\epsilon}})\geq 0$ which means that the latter is zero.

Now let $0< \delta <\frac{1}{2 \max(2T^2,1)}-\beta $. For each
$0<\epsilon <\delta$ there exist $(p_{\epsilon},q_{\epsilon})\in
X\times X$ such that
\begin{eqnarray}
I_{\epsilon}(p_{\epsilon},q_{\epsilon}):&=&
\int_0^T \big[ H_{\epsilon}\big( p_{\epsilon}(t),q_{\epsilon}(t)\big) + H_{\epsilon}^*\big(-\dot{q_{\epsilon}} (t),\dot{p_{\epsilon}}(t)\big) +2\dot{q_{\epsilon}}(t)\cdot p_{\epsilon}(t)\big] \, dt\\
  &&\quad + \psi_2 \big( q_{\epsilon}(T)\big) +\psi_2^* \big(  -p_{\epsilon}(T)\big) +\psi_1 \big( p_{\epsilon}(0)\big) +\psi_1^* \big(
  q_{\epsilon}(0))=0 \nonumber
\end{eqnarray}
We shall show that $(p_{\epsilon},q_{\epsilon})$ is bounded in
$X\times X$. Indeed, similar to the proof of Lemma 3, we get

\begin{eqnarray}
|2\int_0^T\dot{q_{\epsilon}}\cdot p_{\epsilon}\, dt| \leq \max
\left( 2T^2,1\right)\int_0^T
(|\dot{q_{\epsilon}}|^2+|\dot{p_{\epsilon}}|^2)dt +2T
|p_{\epsilon}(0)|^2.
\end{eqnarray}
Combining (12) and (13), we obtain
\begin{eqnarray}
\int_0^T \big[ H_{\epsilon}\big(
p_{\epsilon}(t),q_{\epsilon}(t)\big) +
H_{\epsilon}^*\big(-\dot{q_{\epsilon}}
(t),\dot{p_{\epsilon}}(t)\big) -\max \left(2T^2,1\right)\int_0^T
(|\dot{q_{\epsilon}}|^2+|\dot{p_{\epsilon}}|^2)dt+ \psi_2 \big( q_{\epsilon}(T)\big) \\
  +\psi_2^* \big(  -p_{\epsilon}(T)\big) +\psi_1 \big( p_{\epsilon}(0)\big) +\psi_1^* \big(
  q_{\epsilon}(0))-2T |p_{\epsilon}(0)|^2\leq 0.
\end{eqnarray}
This inequality and the fact that $H$ and $\psi_i^*,  i=1,2$ are
bounded from below, guarantee the existence of a constant  $C>0$
independent of $ \epsilon$ such that
\begin{eqnarray*}
\left(\frac{1}{2(\beta+\delta)}-\max ( 2T^2,1) \right)\int_0^T
(|\dot{q_{\epsilon}}|^2+ |\dot{p_{\epsilon}}|^2)\, dt
 + \psi_2\big( q_{\epsilon}(T)\big)+\psi_1
\big( p_{\epsilon}(0)\big) -2T |p_{\epsilon}(0)|^2\leq C.
\end{eqnarray*}
The coercivity of $\psi_1$ and $\psi_2$ together with the fact that
$\frac{1}{\beta+\delta}-2\max ( 2T^2,1)>0 $ then implies the boundedness
of $(p_{\epsilon},q_{\epsilon})$  in $X\times X$. Therefore, up to a
subsequence $(p_{\epsilon},q_{\epsilon})\rightharpoonup ({\bar p},
{\bar q})$ in $X\times X$.

 Now we show that
\begin{eqnarray}
\label{passtolimit}
I({\bar p},{\bar q})\leq \liminf_{\epsilon\rightarrow
0}I_{\epsilon}(p_{\epsilon},q_{\epsilon})= 0.
\end{eqnarray}
Indeed,  first note that 
\begin{eqnarray*}
\int_{0}^{T}H^*_\epsilon (\dot {p_{\epsilon}},\dot
{q_{\epsilon}})dt:=\inf\limits_{u,v \in  L^2(0,T;\R
^N)}\int_{0}^{T}\left [ H^*(u,v)+\frac{\| \dot
{p_{\epsilon}}-u\|^2}{2 \epsilon}+\frac{\|\dot {q_{\epsilon}}
-v\|^2}{2 \epsilon}\right]dt,
\end{eqnarray*}
and since $H^*$ is convex and lower semi continues, there exists
$u_{\epsilon},v_{\epsilon} \in  L^2(0,T;\R ^N)$  such that this
infimum attains  at $(u_{\epsilon},v_{\epsilon}),$  i.e.
\begin{eqnarray*}
H^*_\epsilon (\dot {p_{\epsilon}},\dot
{q_{\epsilon}})=\int_{0}^{T}\left [
H^*(u_{\epsilon},v_{\epsilon})+\frac{\| \dot
{p_{\epsilon}}-u_{\epsilon}\|^2}{2 \epsilon}+\frac{\|\dot
{q_{\epsilon}} -v_{\epsilon}\|^2}{2 \epsilon}\right]dt.
\end{eqnarray*}
It follows from (14) and the boundedness of
$(p_{\epsilon},q_{\epsilon})$ in $X\times X,$  that there exists a
constant $C>0$ not dependent on $\epsilon$ such that
\begin{eqnarray*}
H^*_\epsilon (\dot {p_{\epsilon}},\dot
{q_{\epsilon}})=\int_{0}^{T}\left [
H^*(u_{\epsilon},v_{\epsilon})+\frac{\| \dot
{p_{\epsilon}}-u_{\epsilon}\|^2}{2 \epsilon}+\frac{\|\dot
{q_{\epsilon}} -v_{\epsilon}\|^2}{2 \epsilon}\right]dt<C.
\end{eqnarray*}
Since $H^*$ is bounded from below, we have
$\int_{0}^{T}\left [\| \dot {p_{\epsilon}}-u_{\epsilon}\|^2+\|\dot
{q_{\epsilon}} -v_{\epsilon}\|^2\right]dt<4C \epsilon$
which means that $(u_{\epsilon}, v_{\epsilon})\rightharpoonup ({\bar{
\dot p}}, {\bar { \dot q}})$ in $L^2\times L^2.$ Hence
\begin{eqnarray*}
\int_{0}^{T}H^*(\dot {{\bar p}},\dot {{\bar q}})dt\leq
\inf\limits_{\epsilon \rightarrow 0}
\int_{0}^{T}H^*(u_{\epsilon},v_{\epsilon}) dt& \leq &
\inf\limits_{\epsilon \rightarrow 0} \int_{0}^{T} \left [
H^*(u_{\epsilon},v_{\epsilon})+\frac{\| \dot
{p_{\epsilon}}-u_{\epsilon}\|^2}{2 \epsilon}+\frac{\|\dot
{q_{\epsilon}} -v_{\epsilon}\|^2}{2 \epsilon}\right]dt\\
 &=&\inf\limits_{\epsilon \rightarrow 0}\int_{0}^{T}H^*_\epsilon (\dot
{p_{\epsilon}},\dot {q_{\epsilon}})dt.
\end{eqnarray*}
We also have
\begin{eqnarray*}
\int_{0}^{T}H( {{\bar p}},  {{\bar q}})dt \leq \inf\limits_{\epsilon
\rightarrow 0}\int_0^T  H\big( p_{\epsilon}(t),q_{\epsilon}(t)\big)
dt & \leq & \inf\limits_{\epsilon \rightarrow 0}\int_0^T \big[
H\big(
p_{\epsilon}(t),q_{\epsilon}(t)\big)+\epsilon(|p_{\epsilon}(t)|^2+|q_{\epsilon}(t)|^2)\big]
dt \\&=&\inf\limits_{\epsilon \rightarrow 0} \int_0^T
H_{\epsilon}\big( p_{\epsilon}(t),q_{\epsilon}(t)\big)dt.
\end{eqnarray*}
Moreover, $\dot {q_{\epsilon}} \rightharpoonup  \dot {\bar q}$
weakly and $ p_{\epsilon} \rightarrow  {\bar p}$  strongly in $L^2$, thus
$\lim_{\epsilon \rightarrow 0}\int_0^T\dot{q_{\epsilon}}\cdot
p_{\epsilon}\, dt= \int_{0}^{T} \dot{\bar q} \cdot  {\bar p} dt.$
Therefore
\begin{eqnarray*}
I({\bar p},{\bar q})\leq \liminf_{\epsilon\rightarrow 0}
\int_0^T \big[ H_{\epsilon}\big( p_{\epsilon}(t),q_{\epsilon}(t)\big) + H_{\epsilon}^*\big(-\dot{q_{\epsilon}} (t),\dot{p_{\epsilon}}(t)\big)
 +2\dot{q_{\epsilon}}(t)\cdot p_{\epsilon}(t)\big] \, dt\\
 + \psi_2 \big( q_{\epsilon}(T)\big) +\psi_2^* \big(  -p_{\epsilon}(T)\big) +\psi_1 \big( p_{\epsilon}(0)\big) +\psi_1^* \big(
  q_{\epsilon}(0))=0.
\end{eqnarray*}
Since by Lemma (1), $I({\bar p},{ \bar q}
)\geq 0$,  the latter is therefore zero, and it follows that
\begin{eqnarray*}
0&=&I({\bar p},{\bar q})\\
&=&
\int_0^T \big[ H\big( {\bar p}(t),{\bar q}(t)\big) + H^*\big( -\dot{{\bar q}} (t),\dot{{\bar p}}(t)\big) +2\dot{{\bar q}}(t)\cdot {\bar p}(t)\big] \, dt\\
  &&\quad + \psi_2 \big( {\bar q}(T)\big) +\psi_2^* \big(  -{\bar p}(T)\big) +\psi_1 \big( {\bar p}(0)\big) +\psi_1^* \big( {\bar q}(0))\\
  &=&
\int_0^T \big[ H\big( {\bar p}(t),{\bar q}(t)\big) + H^*\big( -\dot{{\bar q}} (t),\dot{{\bar p}}(t)\big) +\dot{{\bar q}}(t)\cdot {\bar p}(t)-\dot{{\bar p}}(t)\cdot {\bar q}(t)\big] \, dt\\
  &&\quad + \big[\psi_2 \big( {\bar q}(T)\big) +\psi_2^* \big(  -{\bar p}(T)\big)+{\bar p}(T)\cdot {\bar q}(T)
  \big]\\
  &&\quad + \big[\psi_1 \big( {\bar p}(0)\big) +\psi_1^* \big( {\bar q}(0)\big) -{\bar p}(0)\cdot {\bar q}(0) \big].
  \end{eqnarray*}
The result is now obtained  from the following  3 identities and from the
limiting case in Legendre-Fenchel duality:
\[
H\big( {\bar p}(t),{\bar q}(t)\big) + H^*\big( -\dot{{\bar q}}
(t),\dot{{\bar p}}(t)\big) +\dot{{\bar q}}(t)\cdot {\bar
p}(t)-\dot{{\bar p}}(t)\cdot {\bar q}(t) =0,
\]
\[
\psi_2 \big( {\bar q}(T)\big) +\psi_2^* \big(  -{\bar
p}(T)\big)+{\bar p}(T)\cdot {\bar q}(T)=0,
  \]
  \[
 \psi_1 \big( {\bar p}(0)\big) +\psi_1^* \big( {\bar q}(0)\big) -{\bar p}(0)\cdot {\bar q}(0)=0.
 \]
 \section{The Cauchy problem for Hamiltonian systems}

Here is our result for the corresponding Cauchy problem.
\begin{theorem}
Suppose $H:\R ^N\times\R ^N\ra\R$ is a proper convex lower
semi-continuous function such that
$H(p,q)\rightarrow\infty$ as $|p|+|q|\rightarrow\infty$. Assume that
 \begin{eqnarray}
-\alpha \leq H(p,q)\leq \beta (|p|^r+ |q|^r+1) \quad \quad
(1<r<\infty)
\end{eqnarray}
 where  $\alpha, \beta$ are positive constants.
 Then
the infimum of the  functional
\begin{eqnarray*}
J(p,q):=\int_0^T\left[ H\big( p(t),q(t)\big) +H^*\big(
-\dot{q}(t),\dot{p}(t)\big) +\dot{q}(t)\cdot p(t)-\dot{p}(t)\cdot
q(t)\right]\, dt
\end{eqnarray*}
on the set $D:=\{ (p,q)\in X\times X; \, p(0)=p_0,q(0)=q_0\}$ is equal to zero and is
attained at a solution of
\begin{eqnarray}
\label{Cauchy}
\left\{\begin{array}{rcll}
\dot{p}(t) &\in &\partial_2 H\big( p(t),q(t)\big)\quad &t\in (0,T),\\
-\dot{q} (t) &\in &\partial_1 H\big( p(t),q(t)\big) &t\in (0,T)\\
\hfill (p(0), q(0))&=&(p_0, q_0).
\end{array}\right.
\end{eqnarray}
\end{theorem}
To prove Theorem 2, we first consider the subquadratic case ($1<r<2$).
\begin{proposition}
Assume $H$ is a proper convex and lower semi continuous Hamiltonian that is subquadratic on $\R ^N\times\R ^N$.
Then the infimum of the  functional
\begin{eqnarray}
J(p,q):=\int_0^T\left[ H\big( p(t),q(t)\big) +H^*\big(
-\dot{q}(t),\dot{p}(t)\big) +\dot{q}(t)\cdot p(t)-\dot{p}(t)\cdot
q(t)\right]\, dt
\end{eqnarray}
on $D:=\{ (p,q)\in X\times X,p(0)=p_0,q(0)=q_0\}$ is zero and is
attained at a solution of  (\ref{Cauchy}).
\end{proposition}

\noindent {\bf Proof of Proposition 1:} As in the proof of Theorem 1, it is clear that $J(p,q) \geq 0$ for every $(p,q)\in X\times X$. For the reverse inequality, we may consider --as in section 1-- a perturbed Hamiltonian $H_{\epsilon}$ to insure coercivity, and  then pass to a limit when $\epsilon \to 0$. We therefore can and shall assume that $H$ is coercive.
We then introduce the following
Hamiltonian
\begin{eqnarray*}
L(r,s;p,q):=\int_0^T\left[ (\dot{r},-\dot{s})\cdot
(q,p)+H^*(-\dot{q},\dot{p})-H^*(-\dot{s},\dot{r})+\dot{q}(t)\cdot
p(t)-\dot{p}(t)\cdot q(t)\right]\, dt.
\end{eqnarray*}
and we show that $I({p},{q})\leq\sup\limits_{(r,s)\in
D}L(r,s;{p},{q})$.\\
Indeed, setting
\begin{eqnarray*}
A:=\left\{ (r,s)\in D:r(t)=\int_0^t f(\alpha)\, d\alpha+p_0,s(t)=\int_0^t
g(\alpha)\, d\alpha+q_0,\mbox{ for some }f,g\in L_2(0,T;\R ^N)\right\}.
\end{eqnarray*}
we have
\begin{eqnarray*}
\sup_{(r,s)\in D}L(r,s;{p},{q})&\geq&
\sup_{(r,s)\in A}L(r,s;{p},{q})\\ &=& \sup_{f,g\in L^2}\left\{\int_0^T\left[ (f,-g)\cdot ({q},{p})+H^*(-\dot{{q}},\dot{{p}})-H^*(-g,f)+\dot{{q}}{p}-\dot{{p}}\cdot{q}\right]\, dt\right\}\\
&=& \int_0^T\left[ H^*\left( -\dot{{q}},\dot{{p}}\right) +H({p},{q})+\dot{{q}}\cdot{p}-\dot{{p}}\cdot{q}\right]\, dt\\
&=& I({p},{q}).
\end{eqnarray*}
 The rest follows in the same way as in the proof of Theorem 1, that is the subquadraticity of $H$ gives the right coercivity for $L$ and we are able to apply Ky-Fan's min-max principle as in Theorem 1, to find $(\bar p, \bar q)\in D$ such that $J(\bar p, \bar q)=0$. \\

Now, we deal with the general case, that is when (17) holds  
with $r > 2$. For that we shall use an unusual variation of the standard inf-convolution procedure to reduce the problem  to the subquadratic case where Proposition~(1) applies.\\
For every $\la >0$, define
\begin{eqnarray}
H_\la (p,q):=\inf\limits_{u,v\in\R ^N}\left\{ H(u,v)+\frac{\|
p-u\|_s^s}{s \la^s}+\frac{\| q-v\|_s^s}{s \la^s}\right\}
\end{eqnarray}
where $s=\frac{r}{r-1}$. Obviously, $1<s<2$, and  since $H$ is
convex and lower-semi continuous, the infimum in (20) is attained,
so that for every $p,q \in \R ^N$, there exist unique points $i(p),
j(q) \in \R ^N $ such that
\begin{eqnarray}
H_\la (p ,q )=H\big( i (p ),j (q )\big) +\frac{\| p -i (p
)\|_s^s}{s \la^s}+\frac{\| q -j (q )\|_s^s}{s \la^s}.
\end{eqnarray}

\begin{lemma} The regularized Hamiltonian $H_\la$ satisfies the following properties:
\begin{description}
\item{(i)} $H_\la (p,q)\ra H(p,q)$ as $\la\ra 0^+$.
\item{(ii)} $H_\la (p,q)\leq H(0,0)+\frac{\| q\|_s^s+\| p\|_s^s}{s \la^s}$.
\item{(iii)} $H_\la^* (p,q)=H^*(p,q)+\frac{\la ^r}{r}\left(\| p\|_r^r +\| q\|_r^r\right)$.
\end{description}
\end{lemma}
\noindent {\bf Proof:} (i) and (ii) are easy. For (iii), we have
\begin{eqnarray*}
H_\la^* (p,q) &=& \sup_{u,v\in\R ^N} \left\{ u\cdot p+v\cdot q -H_\la (u,v)\right\}\\
&=& \sup_{u,v\in\R ^N}\left\{ u\cdot p+v\cdot q-\inf_{z,w\in\R ^N}\left\{ H(z,w)+\frac{\| z-u\|_s^s+\| w-v\|_s^s}{s \la^s}\right\}\right\}\\
&=& \sup_{u,v\in\R ^N}\sup_{z,w\in\R ^N}\left\{ u\cdot p+v\cdot q-H(z,w)-\frac{\| z-u\|_s^s}{s \la^s}-\frac{\| w-v\|_s^s}{s_{\la_s}}\right\}\\
&=& \sup_{z,w\in\R ^N}\sup_{u,v\in\R ^N}\left\{ (u-z)\cdot p+(v-w)\cdot q+z\cdot p+w\cdot q-H(z,w)-\frac{\| z-u\|_s^s +\| w-v\|_s^s}{s \la^s}\right\}\\
&=& \sup_{z,w\in\R ^N}\sup_{u_1,v_1\in\R ^N}\left\{ u_1\cdot p+v_1\cdot q-\frac{\| u_1\|_s^s}{s \la^s}-\frac{\| v_1\|_s^s}{s \la^s}+z\cdot p+w\cdot q-H(z,w)\right\}\\
&=& \sup_{u_1,v_1\in\R ^N}\left\{ u_1\cdot p+v_1\cdot q-\frac{\| u_1\|_s^s}{s \la^s}-\frac{\| v_1\|_s^s}{s \la^s}\right\} +\sup_{z,w\in\R ^N}\left\{ z\cdot p+w\cdot q-H(z,w)\right\}\\
&=& \frac{\la^r}{r}\left(\| p\|_r^r+\| q\|_r^r\right)
+H^*(p,q).\quad \square
\end{eqnarray*}

Now consider the Cauchy problem associated to $H_\lambda$.
 By Proposition~(1), there exists $(p_\la ,q_\la )\in X\times X$ such
that $p_\la (0)=p_0,q_\la (0)=q_0$ and
\begin{eqnarray}
0=I(p_\la ,q_\la )=\int_0^T\left[ H_\la (p_\la ,q_\la )+H_\la
^*(-\dot{q}_\la ,p_\la )+\dot{q}_\la \cdot p_\la -\dot{p}_\la \cdot
q_\la\right]\, dt
\end{eqnarray}
yielding
\begin{eqnarray*}
&&\left\{\begin{array}{rcl}
\dot{p}_\la &\in &\partial_2 H_\la (p_\la ,q_\la )\\
-\dot{q}_\la &\in &\partial_1 H_\la (p_\la ,q_\la )\end{array}\right.\\
&&p_\la (0)=p_0,q_\la (0)=q_0.\nonumber
\end{eqnarray*}
From (21), we have
\begin{eqnarray}
H_\la (p_\la ,q_\la )=H\big( i_\la (p_\la ),j_\la (q_\la )\big)
+\frac{\| p_\la -i_\la (p_\la )\|_s^s}{s \la^s}+\frac{\| q_\la
-j_\la (q_\la )\|_s^s}{s \la^s}.
\end{eqnarray}
We now relate $(p_\la, q_\la)$ to the original Hamiltonian.

\begin{lemma}For every $\lambda >0$, we have
\begin{eqnarray*}
\left\{\begin{array}{rcl}
\dot{p}_\la &\in &\partial_2 H \big( i_\la (p_\la ),j_\la (q_\la )\big)\\
-\dot{q}_\la &\in &\partial_1 H \big( i_\la (p_\la ),j_\la
(q_\la )\big).\end{array}\right.
\end{eqnarray*}
\end{lemma}
\noindent {\bf Proof:} From (22) and the definition of
Legendre-Fenchel duality we can write
\begin{eqnarray}
H_\la (p_\la ,q_\la )+H_\la ^*(-\dot{q}_\la ,p_\la )+\dot{q}_\la
\cdot p_\la -\dot{p}_\la\cdot q_\la =0\quad \forall t\in (0,T).
\end{eqnarray}
Part (iii) of Lemma~(5), together with (23) and (24), give
\begin{eqnarray}
0&=&H\big( i_\la (p_\la ),j_\la (q_\la )\big) +\frac{\| p_\la -i_\la (q_\la )\|_s^s +\| q_\la -j_\la (q_\la )\|_s^s}{s \la^s}\nonumber\\
& &+ H^*(-\dot{q}_\la ,\dot{p}_\la )+\frac{\la^r}{r}\big(\|\dot{p}_\la\|_r^r +\|\dot{q}_\la\|_r^r\big)\\
&&+ \dot{q}_\la \cdot p_\la -\dot{p}_\la\cdot q_\la.\nonumber
\end{eqnarray}

\noindent Note that
\begin{eqnarray}
p_\la \cdot\dot{q}_\la =\big( p_\la -i_\la (p_\la
)\big)\cdot\dot{q}_\la +i_\la (p_\la )\cdot\dot{q}_\la \quad
 \text{and} \quad \dot{p}_\la\cdot q_\la =\big( q_\la -j_\la
(q_\la)\big)\cdot\dot{p}_\la +\big(\dot{p}_\la\cdot j_\la
(q_\la)\big).
\end{eqnarray}
By Young's inequality, we have

\begin{eqnarray}
\left|\big( p_\la -i_\la (p_\la )\big)\cdot\dot{q}_\la\right| &\leq& \frac{\| p_\la -i_\la (p_\la )\|_s^s}{s \la^s}+\frac{\la^r}{r}\|\dot{q}_\la\|_r^r\\
\left|\big( q_\la -j_\la (q_\la )\big)\cdot\dot{p}_\la\right|
&\leq&\frac{\| q_\la -j_\la (q_\la
)\|_s^s}{s \la^s}+\frac{\la^r}{r}\|\dot{p}_\la\|_r^r
\end{eqnarray}
Combining (25) - (28) gives
\begin{eqnarray*}
0&=&H\big( i_\la (p_\la ),j_\la (q_\la )\big) +\frac{\| p_\la -i_\la (q_\la )\|_s^s +\| q_\la -j_\la (q_\la )\|_s^s}{s \la^s}\\
&&\quad +H^*(-\dot{q}_\la ,\dot{p}_\la )+\frac{\la^r}{r} \left(\|\dot{p}_\la\|_r^r +\|\dot{q}_\la\|_r^r\right)\\
&&\quad +\big( p_\la -i_\la (p_\la )\big)\cdot\dot{q}_\la +i_\la (p_\la )\cdot\dot{q}_\la -(q_\la -j_\la (q_\la )\cdot\dot{p}_\la -\dot{p}_\la\cdot j_\la (q_\la )\\
&\geq&H\big( i_\la (p_\la ),j_\la (q_\la )\big) +\frac{\| p_\la -i_\la (q_\la )\|_s^s +\| q_\la -j_\la (q_\la )\|_s^s}{s \la^s}\\
&&\quad +H^*(-\dot{q}_\la ,\dot{p}_\la )+\frac{\la^r}{r} \left(\|\dot{p}_\la\|_r^r +\|\dot{q}_\la\|_r^r\right)\\
&&\quad + i_\la (p_\la )\cdot\dot{q}_\la -\dot{p}_\la\cdot j_\la (q_\la )-\frac{\la^r}{r}\left(\|\dot{p}_\la\|_r^r +\|\dot{q}_\la\|_r^r\right) -\frac{\| p_\la -i_\la (p_\la )\|_s^s +\| q_\la -j_\la (q_\la )\|_s^s}{s \la^s}\\
&=& H\big( i_\la (p_\la ),j_\la (q_\la )\big) +H^*(-\dot{q}_\la
,\dot{p}_\la )+i_\la (p_\la )\cdot\dot{q}_\la -\dot{p}_\la\cdot
j_\la (q_\la ).
\end{eqnarray*}
On the other hand, by the definition of Fenchel-Legendre duality
\begin{eqnarray}
H\big( i_\la (p_\la ),j_\la (q_\la )\big) +H^*(-\dot{q}_\la
,\dot{p}_\la )+i_\la (p_\la )\cdot\dot{q}_\la -\dot{p}_\la\cdot
j_\la (q_\la )\geq 0
\end{eqnarray}
which means we have equality in (29), so that
\begin{eqnarray*}
\left\{\begin{array}{rcl}
\dot{p}_\la &\in &\partial_2 H\big( i_\la (p_\la ),j_\la (q_\la )\big),\\
-\dot{q}_\la &\in &\partial_1 H\big( i_\la (p_\la ),j_\la (q_\la
)\big).\quad\square\end{array}\right.
\end{eqnarray*}

\begin{lemma} With the above notation we have:
\begin{enumerate}
\item   $ \sup\limits_{t\in (0,T)}\left| q_\la -j_\la (q_\la )\right|
+|p_\la -i_\la (p_\la )|\leq c \la$, where $c$ is a constant.

\item   If $H(p,q)\rightarrow\infty$  as $|p|+|q|\rightarrow\infty$
then $ \sup\limits_{t\in (0,T), \la >0}| q_\la |+|j_\la (q_\la )|
+|p_\la|+|i_\la (p_\la )| < \infty.$
\end{enumerate}
\end{lemma}
 {\bf Proof:}  For every $\lambda >0$ and
$t \in (0,T)$, we have
\begin{eqnarray*}
\left\{\begin{array}{rcl}
\dot{p}_\la &= &\partial_2 H_\la (p_\la ,q_\la )\\
-\dot{q}_\la &= &\partial_1 H_\la (p_\la ,q_\la )\end{array}\right.
\end{eqnarray*}
Multiplying the first equation by $\dot{q}_\la$ and the second one
by $\dot{p}_\la$, give
\begin{eqnarray*}
\left\{\begin{array}{rcl}
\dot{p}_\la\dot{q}_\la &= &\dot{q}_\la\partial_2 H_\la (p_\la ,q_\la )\\
-\dot{q}_\la p_\la &= &\dot{p}_\la\partial_1 H_\la (p_\la ,q_\la
)\end{array}\right.
\end{eqnarray*}
So $\frac{d}{dt}H_\la (p_\la ,q_\la )=0$ and
$H_\la\big( p_\la (t),q_\la (t)\big)=H_\la\big( p(0),q(0)\big)\leq
H\big( p(0),q(0)\big) :=c<+\infty.$
Hence, it follows from (23) that
\begin{eqnarray*}
H\big( i_\la (p_\la (t),j_\la (q_\la (t)\big)+\frac{\| p_\la -i_\la
(p_\la )\|_s^s +\| q_\la -j_\la (q_\la )\|_s^s}{s \la^s}\leq c
\end{eqnarray*}
which yields
$\sup_{t\in (0,T]}|q_\la -j_\la (q_\la )|+|p_\la -i_\la p_\la |\leq
c\la$
and
$\sup_{t\in (0,T]}H\left( i_\la \big( p_\la (t)\big) ,j_\la\big(
q_\la (t)\big)\right) <+\infty.
$
Since $H$ is coercive the last equation gives
$ \sup_{t\in (0,T) \la}|j_\la (q_\la )|+|i_\la p_\la | <\infty,$
which together with Part (i) prove the lemma. \hspace{7in} $\square$

\begin{lemma}  We have the following estimate:
\begin{eqnarray*}
 \sup\limits_{t\in [0, T],\la>0}|\dot{p}_\la(t)|+|\dot{q}_\la (t)|<+\infty.
 \end{eqnarray*}
\end{lemma}
{\bf Proof:} Since
$-\alpha < H(p,q)\leq \beta |p|^r+ \beta|q|^r+ \beta$   where $r>2$,
an easy calculation shows that if $(p^*,q^*) \in  \partial H(p,q)$ then
\begin{equation}
|p^*|+|q^*|\leq \left \{ s (2\beta)^{\frac
{r}{s}}(|p|+|q|+\alpha+\beta)+1 \right \}^{r-1}
\end{equation}
Since by  Lemma (6) we have $(\dot{p}_\la,-\dot{q}_\la) =
\partial H(i_\la (p_\la ),j_\la (q_\la ))$, it follows from (30) that
\begin{equation*}
|\dot{p}_\la|+|\dot{q}_\la|\leq \left \{ s (2\beta)^{\frac
{r}{s}}(|i_\la (p_\la )|+|j_\la (q_\la )|+\alpha+\beta)+1 \right
\}^{r-1}
\end{equation*}
which together with Lemma 7 prove the desired result.\\

\noindent{\bf End of proof of Theorem 2:} From Lemma~(6), we have
\begin{eqnarray}
\int_0^T\left[ H\big( i_\la (p_\la ),j_\la (q_\la ))+H^*(-\dot{q}_\la
,\dot{p}_\la )+\dot{q}_\la\cdot i_\la (p_\la ) -\dot{p}_\la\cdot
j_\la (q_\la )\right]\, dt=0.
\end{eqnarray}
while $p_\la (0)=p_0$ and $q_\la (0)=q_0$. By Lemma~(8),
$\dot{p}_\la$ and $\dot{q}_\la$ are bounded in $L^2(0,T;\R ^N)$ so
there exists $({p},{q})\in X\times X$ such that
$\dot{p}_\la\rightharpoonup\dot{{p}}$ and
$\dot{q}_\la\rightharpoonup\dot{{q}}$ weakly in $L^2(0,T;\R ^N)$ and
$p_\la\ra{p}$ and $q_\la\ra{q}$ strongly in $L_\infty (0,T;\R ^N)$.
So by Lemma~(7), $i_\la (p_\la )\ra{p}$ and $j_\la(q_\la )\ra{q}$
strongly in $L_\infty (0,T;\R ^N)$. Hence by letting $\la\ra 0$ in
(31), we get
\begin{eqnarray*}
\int_0^T\left[ H({p},{q})+H^*\left(
-\dot{{q}},\dot{{p}}\right) +\dot{{q}}\cdot{p}
-\dot{{p}}\cdot{q}\right]\, dt\leq 0,
\end{eqnarray*}
which means ${p}(0)=p_0$ and ${q}(0)=q_0$ and 
\begin{eqnarray*}
\left\{\begin{array}{rcl}
\dot{p}(t) &= &\partial_2 H({p(t)},{q(t)})\\
-\dot{q}(t) &= &\partial_1 H({p(t)},{q(t)}).\end{array}\right.
\end{eqnarray*}

\section{ Semi-Convex Hamiltonian systems}

In this section, we consider the following system:
\begin{eqnarray}
\left\{\begin{array}{rcll}
\dot{p}(t) &\in &\partial_2 H\big( p(t),q(t)\big) +\de _1 q(t)\quad &t\in (0,T)\\
-\dot{q}(t) &\in &\partial_{1} H\big( p(t),q(t)\big) +\de _2 p(t)
&t\in (0,T)\\
q(0)&\in&\partial\psi_1\big( p(0)\big)\\
-p(T)&\in&\partial\psi_2\big( q(T)\big)
\end{array}\right.
\end{eqnarray}
where $\de_1,\de_2 \in\R$. Note that if $\delta_i \geq 0$ then the problem reduces to the one studied in section 1 with a new convex Hamiltonian $\tilde H (p,q)= H(p,q)+\frac{\delta_1}{2}|q|^2 +\frac{\delta_2}{2}|p|^2$.  The case that concerns us here is when $\delta_i <0$.  
\begin{theorem}
Suppose $H:\R ^N\times\R ^N\ra\R$ is a proper convex lower
semi-continuous Hamiltonian that is $\beta$-subquadratic with 
\begin{equation}
\label{betaprime}
\beta < \frac
{1}{4} \min \{\frac{1-4T^2|\de_1|^2}{\max(2T^2,1)-2\de_1 T^2},
\frac{1-4T^2|\de_2|^2}{\max(2T^2,1)-2\de_2 T^2} \}.
\end{equation}
Assume 
\begin{equation}
\hbox{$|\de_i |<\frac{1}{2T}$ for $i=1,2$}, 
\end{equation}
and let $\psi_1$ and $\psi_2$ be  convex lower semi-continuous functions on $\R ^N$
satisfying
\begin{eqnarray}\begin{array}{l}
\liminf\limits_{|p|\ra +\infty}\frac{\psi_1
(p)}{|p|^2}>\frac{T|\de_2|^2}{\beta}+2T(1-\de_2)\quad {\rm and }\quad 
\liminf\limits_{|p|\ra +\infty}\frac{\psi_2
(p)}{|p|^2}>\frac{T|\de_1|^2}{\beta}-2T\de_1.
\end{array}
\end{eqnarray}
Then the minimum of
the functional
\begin{eqnarray*}
I(p,q):& =& \int_0^T \big[
H\big(p(t),q(t)\big)+H^*\big(-\dot{q}(t)-\de_2 p(t),\dot{p}(t)-\de_1
q(t)\big) +\de_1 |q|^2 +\de_2|p|^2+2\dot{q}(t)\cdot p(t)\big] \, dt\\
 & & + \psi_2 \big(  q(T)\big)+\psi_2^*\big(-p(T)\big)
+ \psi_1 \big( p(0)\big) +\psi_1^* \big( q(0)\big)
  \end{eqnarray*}
on $Y=X\times X$ is equal to zero and is attained at a solution of (32).
\end{theorem}
\bigskip
By considering a perturbed Hamiltonian $H_{\epsilon}$, then passing to a limit when $\epsilon \to 0$ as in section 1, we can and shall  assume that $H$ is coercive.
Also,  note that for every $(p,q) \in Y,$
\begin{eqnarray*}
I(p,q)& =& \int_0^T \big[
H\big(p(t),q(t)\big)+H^*\big(-\dot{q}(t)-\de_2 p(t),\dot{p}(t)-\de_1
q(t)\big) +\de_1 |q|^2 +\de_2|p|^2+\dot{q}(t)\cdot p(t)-\dot{p}(t)\cdot q(t)\big] \, dt\\
 & & + \big[\psi_2 \big(  p(T)\big)+\psi_2^*\big(-q(T)\big)+p(T)\cdot q(T)  \big]\\
 & & + \big[\psi_1 \big( p(0)\big) +\psi_1^* \big( q(0)\big)-p(0)\cdot q(0) \big]\\
 &\geq&0,
  \end{eqnarray*}
by three applications of Legendre inequality. \\
For the reverse inequality, we
 introduce the following functional Lagrangian
$L:Y\times Y\ra\R$ defined by
\begin{eqnarray*}
L(r,s;p,q): & = & \int_0^T\Big[ H^* (-\dot{q}-\de_2 p,\dot{p}-\de_1 q) - H^*(-\dot{s}-\de_2 r,\dot{r}-\de_1 s)\\
&&\quad \quad+(\dot{r}-\de_1s,-\dot{s}-\de_2 r)\cdot (q,p)+\de_1 |q|^2+\de_2 |p|^2 +2\dot{q}\cdot p\Big]\, dt\\
& & - p(T)\cdot s(T)+\psi_2 \big( q(T)\big) -\psi_2\big(s (T)\big)
+r(0)\cdot q(0)+\psi_1\big( p(0)\big)-\psi_1\big( r(0)\big).
\end{eqnarray*}
 In order to apply the anti-selfduality argument, we   need the following Lemma.

\begin{lemma}
For any $f,g\in L^2(0,T;\R ^N)$ and $x,y\in\R ^N$, there exists
$(r,s)\in X\times X$ such that
\begin{eqnarray}
\left\{\begin{array}{rcl}
\dot{r}(t) &= &\de_2 s(t)+f(t)\\
-\dot{s}(t) &= &\de_1 r(t)+g(t)\\
r(0) &= &x\\
s(T)&=&y.\end{array}\right.
\end{eqnarray}
\end{lemma}
 {\bf Proof:} This is standard and is essentially a  linear system of ordinary differential equations.
 Also, one can  rewrite the problem as
follows.
\begin{eqnarray*}
\left\{\begin{array}{rcl}
-\dot{r}(t)+f(t) &= &-\partial_2G\big( r(t),s(t)\big)\\
\dot{s}(t)+g(t) &= &-\partial_1G\big( r(t),s(t)\big)\\
r(0) &= &x\\
s(T)&=&y.
\end{array}\right.
\end{eqnarray*}
where $G\big( r(t),s(t)\big) =-\frac{\de_1}{2}\int_0^T |r(t)|^2\,
dt-\frac{\de_2}{2}\int_0^T |s(t)|^2\, dt$. Hence
\begin{eqnarray*}
G^*\big(\dot{s}(t)+g(t),-\dot{r}(t)+f(t)\big) =-\frac{1}{2\de_2}
   \int |\dot{s}(t)+g(t)|^2 -\frac{1}{2\de_1}\int |-\dot{r}(t)+f(t)|^2\, dt
\end{eqnarray*}
One can show as in Theorem 1 that whenever $|\de_i |<\frac{1}{2T}$,
coercivity holds and the following infimum is achieved at a solution
of (36).
\begin{eqnarray*}
0&=&\inf\limits_{(r,s)\in D\subseteq X\times X}
  G^*\big(\dot{s}(t)+g(t),-\dot{r} (t)+f(t)\big)
  +G\big(r(t),s(t)\big)\\
   &&\quad \quad \quad \quad +\int_0^T \dot{r}(t)\cdot s(t)\, dt
  -\int_0^T \dot{s}(t)\cdot r(t)\, dt\\
  &&\quad \quad \quad \quad -\int_0^T \big( f(t)\cdot s(t)+r(t)\cdot g(t)\big)\, dt
\end{eqnarray*}
where $D=\left\{ (r,s)\in X\times X\mid r(0)=x,s(T)=y\right\}$. \\

\begin{lemma} For every $({p},{q})\in
X\times X$, we have
\begin{eqnarray*}
I({p},{q})\leq \sup\limits_{(r,s)\in X\times X}
L(r,s;{p},{q}).
\end{eqnarray*}
\end{lemma}
{\bf Proof:}  
 Use the above lemma to write
 \begin{eqnarray*}
\sup\limits_{(r,s)\in X\times X}L(r,s;{p},{q})
  &= & \sup\limits_{f,g\in L^2}\sup\limits_{x,y\in\R ^N}
  \int_0^T\Big[ (f,g)\cdot ({q},{p})+H^*(-\dot{{q}}
  -\de_2 {p},\dot{{p}} -\de_1{q} )\\
& & \quad \quad  \quad \quad -H^*(g,f)+\de_1 |q|^2+\de_2 |p|^2+2\dot{q}\cdot p\Big]\, dt\\
& & - {p} (T)\cdot y+\psi_2\big({q} (T)\big)
    - \psi_2 (y)+x\cdot{q} (0)+\psi_1\big({p} (0)\big) -\psi_1(x)\\
&=& \int_0^T\left[ H^*\left( -\dot{{q}}-\de_2
{p},\dot{{p}}-\de_1{q}\right)
   +H({q},{p})+\de_1 |q|^2+\de_2 |p|^2+2\dot{q}\cdot p\right]\, dt\\
& & + \psi_2\big({q} (T)\big) +\psi_2^* \big( -{p}(T)\big)
    +\psi_1\big({p} (0)\big) +\psi_1^*\big({q}
    (0)\big)\\
    &=&I({p},{q}).
\end{eqnarray*}
In order to again apply Ky-Fan's lemma, it remains to establish the following coercivity property.

\begin{lemma} Under the above hypothesis, we have
\begin{equation*}
\hbox{
$L(0,0;p,q)\ra +\infty$ as $ \| p\| +\| q\| \ra +\infty.$}
\end{equation*}
\end{lemma}
  {\bf Proof:} Since $H, \psi^*_1, \psi^*_2$ are bounded from below and  $ \frac {1}{2(\beta+\epsilon)}(|p|^2+|q|^2)-\gamma \leq
H^*_{\epsilon}(p,q)$, modulo a constant we have,
 \begin{eqnarray}
L(0,0;p,q)\geq \frac {1}{2(\beta+\epsilon)}\int_0^T \left(
|\dot{q}+\de_2 p|^2 +|\dot{p}-\de_1 q|^2\right)\, dt
& + & \int_0^T \left(\de_1 |q|^2 +\de_2 |p|^2+2\dot{q}\cdot p\right)\, dt  \nonumber\\
& + & \psi_2 \big( q(T)\big) +\psi_1\big( p(0)\big).
\end{eqnarray}
It is easily seen that
\begin{eqnarray}
| \dot{q}+\de_2 p|^2\geq \frac{1}{2} |\dot{q}|^2 -|\de_2|^2 | p|^2,
\quad \text{and} \quad | \dot{p}-\de_1 q|^2\geq \frac{1}{2}
|\dot{p}|^2 -|\de_1|^2 | q|^2.
\end{eqnarray}
It follows from (9) that,
\begin{eqnarray}
\int_{0}^{T}| p(t)|^2 dt\leq 2\left( T^2\int_{0}^{T}| \dot{p}|^2dt
+T|p(0)|^2 \right) \quad \text{and} \quad \int_{0}^{T} | q(t)|^2dt
\leq 2\left( T^2\int_{0}^{T} |\dot{q}|^2dt  +T|q(T)|^2 \right) .
\end{eqnarray}
Combining (37) and  (38) gives,
 \begin{eqnarray}
\int_{0}^{T}\left[| \dot{q}+\de_2 p|^2+| \dot{p}-\de_1 q|^2\right]
dt
  &\geq & \int_{0}^{T}\left[\frac{1}{2}\left( |\dot{q}|^2 +|\dot{p}|^2\right) -|\de_2|^2 | p|^2 -|\de_1|^2 | q|^2 \right]dt\nonumber \\
  &\geq &\int_{0}^{T} \left[\frac{1}{2}\left( |\dot{q}|^2 +|\dot{p}|^2\right) -2T^2\big(|\de_2|^2 |\dot{ p}|^2 +|\de_1|^2 | \dot{q}|^2\big)\right] dt\\
  \quad && -2T \big( |\de_2|^2|p(0)|^2+|\de_1|^2|q(T)|^2\big) \nonumber  \\
  &\geq & \int_{0}^{T} \frac{1}{2} \left[ (1-4T^2|\de_2|^2) |\dot{p}|^2 +(1-4T^2|\de_1|^2 )|\dot{q}|^2\right]dt\\
  \quad && -2T \big( |\de_2|^2|p(0)|^2+|\de_1|^2|q(T)|^2\big)   \nonumber\\
  & = & \int_{0}^{T}\frac{1}{2}\left[ \epsilon_1|\dot{q}|^2 +\epsilon_2|\dot{p}|^2\right]-2T \big( |\de_2|^2|p(0)|^2+|\de_1|^2|q(T)|^2\big)
\end{eqnarray}
 where $\epsilon_i:=1-4T^2|\de_i|^2>0$ since   
$|\de_i |<\frac{1}{2T}$.\\
Also, similarly to  the proof of Theorem~(1), we
get
\begin{eqnarray}
\left| \int_0^T 2\dot{q}\cdot p\, dt\right| \leq \max  \left(
2T^2,1\right) \int_0^T\left( |\dot{q}|^2+|\dot{p}|^2\right)\, dt
+2T|p(0)|^2.
\end{eqnarray}
Hence, combining (37)-(41)  yields
\begin{eqnarray*}
L(0,0;p,q) & \geq & \frac {1}{4(\beta+\epsilon)}\int_0^T
   \left[ \epsilon_1 |\dot{q}|^2 +\epsilon_2 |\dot{p}|^2\right]
   -A(\de_1,T)\int_0^T  |\dot{q}|^2\, dt-A(\de_2,T)\int_0^T  |\dot{p}|^2\, dt\\
   &  & + \psi_2 \big( q(T)\big) -\frac {T|\de_1|^2}{\beta+\epsilon}|q(T)|^2+2\de_1T |q(T)|^2\\
   &  & + \psi_1\big( p(0)\big) -\frac {T|\de_2|^2}{\beta+\epsilon}|p(0)|^s-2T(1-\de_2) |p(0)|^2
\end{eqnarray*}
where $A(\de_i,T)=\max(2T^2,1)-2\de_i T^2. $ This inequality
together with the coercivity condition on $\psi_1$ and $\psi_2$, and
the fact that $\beta < \frac {1}{4}\min \{\frac{\epsilon_1}{
A(\de_1,T) },\frac{\epsilon_2}{ A(\de_2,T) }   \} $ yield the
claimed result.

\end{document}